\documentclass[11pt]{amsart}
\usepackage{amsmath}
\usepackage{amssymb}
\usepackage[a4paper]{geometry}
\theoremstyle{plain}
\newtheorem{teo}{Theorem}
\newtheorem{pro}{Proposition}
\newtheorem{lem}{Lemma}

\newtheorem{cor}{Corollary}

\theoremstyle{definition}
\newtheorem{defi}{Definition}

\theoremstyle{remark}
\newtheorem*{nota}{Remark}
\newtheorem*{note}{Remarks}

\newcommand{\re}{\mathbb{R}}

\usepackage{mathrsfs}
\title{Free Jacobi Process} 
\date{\today  }
\begin{document}
\maketitle
\centerline{N. DEMNI\footnote{Laboratoire de Probabilit\'es et Mod\`eles Al\'eatoires, Universit\'e de Paris VI, 4 Place Jussieu, Case 188, F-75252 Paris Cedex 05, e-mail: demni@ccr.jussieu.fr.}} 
\begin{abstract}In this paper, we define and study  free Jacobi processes of parameters $\lambda > 0$ and $0< \theta \leq 1$, as the limit of the complex version of the matrix Jacobi process already defined by Y. Doumerc. In the first part, we focus on the stationary case  for which we compute the law (that does not depend on time) and derive, for $\lambda \in ]0,1]$ and $1/\theta \geq \lambda +1$ a free SDE analogous to the classical one. In the second part, we generalize this result under an additional condition. To proceed, we set a  recurrence formula for the moments of the process using free stochastic calculus. This will also be used to compute the p. d. e. satisfied by the Cauchy transform of the free Jacobi's law.    
\end{abstract}
\section{Introduction}
Free stochastic calculus with respect to the free additive Brownian motion was introduced by Kummerer and Speicher and then was developed by Biane and Speicher (cf \cite{Bia1}) who gave a ''free It\^o formula'' for polynomials and a special class of functions. They first established results for simple (piecewise constant) '' bi-processes'' (the prefixe '' bi '' is due to the non-commutativity, namely the integrand can be multiplied both on the right and on the left) belonging to some norm vector space. Next, they extend continuously the results by controlling the operator norm of the stochastic integral, more precisely, they established the $L^{\infty}$ version of the free BDG inequality which extends the $L^p$ free BDG inequality already set by Pisier and Xu (cf \cite{PX}) to the case $p = \infty$. Few years later, M. Capitaine and C. Donati defined and studied  free Wishart processes. A free SDE similar to the classical one (i. e. of the complex Wishart processes) is provided. 
Here, we will define the free Jacobi process $(J_t)_{t \geq 0}$ with parameters $(\lambda, \theta)$ as the limit of the complex matrix Jacobi process with parameters $(p,q)$ when the size of the matrix  goes to infinity. Then, we will derive, for suitable parameters, the free SDE verified by this process. The rest of the paper consists of 5 sections which are devoted to the following topics: in section 2, we introduce the free Jacobi process with parameters $\theta \in ]0,1]$ and $\lambda > 0$, denoted $FJP(\lambda,\theta)$. In section 3, we use the free stochastic calculus and the polar decomposition to derive a free SDE provided that $J$ and $P-J$ remain injective. Once deriving this SDE, we will give conditions ensuring the required injectivity. Section 4 considers the stationary case: we will compute the  $J_t^{'}$s law using Nica and Speicher's result on the compression by a free projection (cf \cite{Nica} or \cite{Spei} p. 87).  Section 5 deals with a more general setting and is concerned with the free Jacobi process $(J_t)_{t \geq 0}$ to which we extend our result assuming in addition that $J_0$ is invertible : to do this, we set a recurrence formula for the Jacobi moments which enables us to get a lower bound of $\tilde{\Phi}(\log(P-J_t))$ for all $t > 0$ where $P$ and $\tilde{\Phi}$ are  respectively the unit and the state of the so-called ''compressed space''.  Note that this free SDE involves a complex free Brownian motion which is the (normalized) limit of the complex non self-adjoint  Brownian motion involved in its classical counterpart. In the last section, we compute the p. d. e. satisfied by the Cauchy transform of the free Jacobi's law. Hopelessly, no results concerning the Jacobi's law has been deduced.    
\section{Complex Matrix Jacobi process and free Jacobi process}
We refer to \cite{Yan} for facts about real matrix Jacobi processes. In the remainder of this paper, we are interested in the complex case. Let $\Theta$ be a $d \times d$ unitary Brownian matrix, i. e, a unitary matrix-valued process such that :\begin{itemize}
\item $\Theta_0 = {\it I}_d$ \\
\item $(\Theta_{t_{i}}\Theta_{t_{i-1}}^{-1} , 1 \leq i \leq n)$ are independent for any collection $0 < t_1< \hdots < t_n$, \\
\item $\Theta_t\Theta_s^{-1}, \, s < t$ have the same distribution as $\Theta_{t-s}$ (cf \cite{Bia}),\end{itemize} 
and $X = \Pi_{m, p}\Theta$ is the left corner of $\Theta$, then $J := XX^{\star}$ is a $m \times m$ complex matrix Jacobi process starting at ${\it I}_m$ and such that $0 \leq X_tX_t^{\star} \leq {\it I_m}$. If $X_t$ is the $m \times p$ left corner of $\tilde{Z}\Theta_t$ where $\tilde{Z}$ is a $d \times d$ unitary random matrix independent of $\Theta$, then $XX^{\star}$ is a $m \times m$ complex matrix Jacobi process starting from $X_0X_0^{\star}$ where $X_0$ is the left corner of $\tilde{Z}$. Let us take a $d(m) \times d(m)$ unitary Brownian matrix $Y_t(d(m))$ and let us consider the left corner $X_t = \Pi_{m,p(m)}Y_t(d(m))$ such that :
\begin{equation*}\lim_{m \rightarrow \infty}\frac{m}{d(m)} = \lambda\theta, \qquad \lim_{m \rightarrow \infty}\frac{m}{p(m)} = \lambda > 0\end{equation*}
Then, we may write:  \begin{equation*} XX_t^{\star} \oplus 0_{d(m) - m} = P_m  Y_t(d(m)) Q_{m} Y_t^{\star}(d(m)) P_m \end{equation*} where $P_m, Q_m$ are two projections defined by: 
\begin{equation*} P_m = \left(\begin{array}{cc }
    {\it I}_m  &    \\
      &   0
\end{array}\right) \qquad Q_m = \left(\begin{array}{ cc}
  {\it I}_{p(m)}    &    \\
      &   0
\end{array}\right)\end{equation*} 
Recall that the unitary Brownian matrix, considered as a non-commutative variable in the non-commutative space $\mathscr{A}_d := \bigcap_{p > 0}L^p\left(\Omega, \mathscr{F}, (\mathscr{F}_t), \mathbb{P}\right) \otimes  \mathscr{M}_d(\mathbb{C})$ equipped with the normalized trace expectation $ \mathbb{E} \otimes tr_d$ , converges in distribution to the free multiplicative Brownian motion (cf \cite{Bia}) $Y$ in some non commutatif space $(\mathscr{A}, \Phi)$, that is : \begin{equation*}
\lim_{d \rightarrow \infty} \mathbb{E}[tr_d (Y_{t_1}(d)\dots Y_{t_n}(d))]  = \Phi(Y_{t_1}\dots Y_{t_n}) \end{equation*}for any collection of times $ t_1 ,\dots, t_n$, which implies that
\begin{equation*}
\lim_{n \rightarrow \infty} \mathbb{E}[tr_d (Y_t(d)^k)]  = \Phi(Y_t^k),\qquad  k \geq 0 \end{equation*}
Recall also that $Y$ is unitary, $Y_0 = {\it I}$ ( the unit of $\mathscr{A}$), has free left increments, that is, for a collection of times $0 < t_1 < t_2 < \hdots < t_n$, $Y_{t_n}Y_{t_{n-1}}^{-1}, \hdots, Y_{t_2}Y_{t_1}^{-1}$ are free and that the law $\nu_t$ of $Y_t$ is given by its $\Sigma$-transform:
\begin{equation*}\Sigma_{\nu_t}(z) = e^{\frac{t}{2}\frac{1+z}{1-z}}, \qquad \nu_{t+s} = \nu_t \boxtimes \nu_s\end{equation*} where $\boxtimes$ denotes the free multiplicative convolution (cf \cite{Ber}, \cite{Bia}). Furthermore, we have : \begin{eqnarray*}\lim_{m \rightarrow \infty}tr_{d(m)}(P_m) &=& \lambda \theta \\
 \lim_{m \rightarrow \infty}tr_{d(m)}(Q_{m}) &=&  \theta, \quad \theta \in ] 0, 1[\\
\lim_{m \rightarrow \infty}\mathbb{E}(tr_m(X_tX_t^{\star})) &=&\lim_{m \rightarrow \infty}\left[\frac{d(m)}{m}\mathbb{E}(tr_{d(m)}(P_mY_t(d(m))Q_mY_t^{\star}(d(m))P_m))\right].
 \end{eqnarray*} 
Furthermore,  if $(U_t(n))_{ t \geq 0}$ is a family of independent unitary random matrices such that the distribution of $U_t(n)$ is equal to that of $VU_t(n)V^{\star}$ for any unitary matrix $V$ (unitary invariant) and such that $U_t(n)$ converges in distribution, and if $(D_t(n))_{t \geq 0}$ is a family of constant matrices converging in distribution and such that $\sup_n||D_t(n)|| < \infty$, then the families  \begin{equation*}
 \{U_s(n), U_s^{\star}(n)\}_s , \{D_t(n), D_t^{\star}(n), t \geq 0\}\end{equation*} are asymptotically free as $n \rightarrow \infty$ (cf \cite{Hiai} Theorem. 4. 3. 1. p. 157). \\ Thus, applying these results to the family $\{Y_t(d(m))Y_s^{-1}(d(m))\}_{0\leq s < t}$ which is unitary invariant since $(Y_t)$ is right-left invariant,  the family $\{P_m, Q_m\}$ and using increments freeness of $Y$ mentionned above, one can see that:
\begin{equation*} \{(Y_t(d(m)))_{t \geq 0}, (Y^{\star}_t(d(m)))_{t \geq 0}\}, \{P_m, Q_{m}\} \end{equation*} are asymptotically free, i. e, its limiting distribution  in $(\mathscr{A}_{d(m)},  \mathbb{E} \otimes tr_{d(m)})$ as $m$ goes to infinity is the distribution of $\{(Y_t)_{t \geq 0}, (Y^{\star}_t)_{t \geq 0}, P, Q\}$ in $(\mathscr{A}, \Phi)$ such that \\
\begin{enumerate}
  \item $Y$ is a free multiplicative Brownian motion in $(\mathscr{A}, \Phi)$.
  \item $P$ is a projection with $\Phi(P) = \lambda \theta \leq 1$, \qquad $\theta \in ]0, 1[$.
  \item $Q$ is a projection with $\Phi(Q) = \theta$.
  \item $QP = PQ = \left\{\begin{array}{ccc}P & \textrm{if} \quad \lambda \leq 1 \\  Q & \textrm{if} \quad \lambda > 1\end{array}\right.$
  \item $\{(Y_t)_{t \geq 0}, (Y_t^{\star})_{t \geq 0}\}$ and $\{P, Q\}$ are free.
\end{enumerate} Hence, we deduce that the limiting distribution of the complex matrix Jacobi process $(J_t)_{t \geq 0}$ in $(\mathscr{A}_m, \mathbb{E} \otimes tr_m)$ is the distribution of $(PY_tQY^{\star}_tP)_{t \geq 0}$ in the compressed non-commutative space $P\mathscr{A}P$ equipped with the state \begin{equation*}
\widetilde{\Phi} = \frac{1}{\Phi(P)} \Phi_{|P\mathscr{A}P} = \frac{1}{\lambda \theta} \Phi_{|P\mathscr{A}P}.\end{equation*}  
Hence, we shall define the free Jacobi process as follows:

\begin{defi}Let $(\mathscr{A}, \Phi)$ be a non-commutative probability space. Let $\theta \in ] 0, 1[$ and $\lambda > 0$ such that $\lambda \theta \leq 1$. Let $P$ and $Q$ be two projections such that \begin{equation*} \Phi(Q) = \theta ,\quad \Phi(P) = \lambda \theta, \textrm{and} \, PQ = QP = \left\{\begin{array}{ccc}P & \textrm{if} \quad \lambda \leq 1  \\  Q & \textrm{if} \quad \lambda > 1\end{array}\right.\end{equation*}
Let $Y$ be a free multiplicative Brownian motion  such that $\{(Y_t)_{t \geq 0}, (Y_t^{\star})_{t \geq 0}\}$ and $\{P, Q\}$ is a free family in $(\mathscr{A}, \Phi)$. We will say that a process $J$ in a non-commutative space $(B, \Psi)$ is a free Jacobi process with parameters $(\lambda, \theta)$  if its distribution in $(B, \Psi)$ is equal to the distribution of the process $(PY_tQY^{\star}_t P)_{t \geq 0}$ in $\displaystyle(P\mathscr{A}P, (1/\Phi(P))\Phi_{|P\mathscr{A}P})$. This process starts at $J_0 = P$ if $\lambda \leq 1$ and $J_0 = Q$ if $\lambda > 1$. 
\\ Equivalently, the law of $J$ is the limiting distribution of a  complex matrix Jacobi process when $\displaystyle \frac{m}{p(m)} \underset{m \rightarrow \infty }{\longrightarrow} \lambda$ and $\displaystyle \frac{m}{d(m)} \underset{m \rightarrow \infty}{\longrightarrow} \lambda \theta$. \end{defi}
We also define the free Jacobi process starting at $J_0$: 
\begin{defi}Let $Y$ be a free multiplicative Brownian motion and $Z$ a unitary operator free with $Y$. Then, the process defined by $\tilde{Y} = YZ$ is a free multiplicative Brownian motion starting at $\tilde{Y_0} = Z$. Moreover, if $Z$ is free with $\{P,Q\}$, then the process $\tilde{J}$ defined by : \begin{equation*} \tilde{J}_t := P\tilde{Y}_tQ\tilde{Y}_t^{\star}P \end{equation*} is called a free Jacobi process  with parameters $(\lambda, \theta)$ and starting at $J_0 = PZQZ^{\star}P$. \end{defi} 
For the sake of simplicity, we will write $Y$ for a free multiplicative Brownian motion starting at $Y_0$ and $J$ for a free Jacobi process starting at $J_0$. 

\section{Free Jacobi Process And Free Stochastic Calculus} 
We refer to \cite{Bia} and \cite{Bia1} for free stochastic calculus and notations. Let $(\mathscr{A}, (\mathscr{A}_t), \Phi)$ be a filtered $W^{\star}$ non-commutative probability space, i. e, 
$\mathscr{A}$ is a von Neumann Algebra, $\Phi$ is a faithful normal tracial state and $(\mathscr{A}_t)$ is an increasing family of unital, weakly closed $\star$-subalgebras of $\mathscr{A}$ .  Recall that a bi-process $(U_t = \sum_i A_t^i \otimes B_t^i) \in \mathscr{A} \otimes \mathscr{A}^{op}$ is adapted if $U_t \in \mathscr{A}_t \otimes \mathscr{A}_t$ for all $t \geq 0$. Furthermore, adapted bi-processes form a complex vector space that we endow with the norm: 
\begin{equation*}
||U|| = \left(\int_0^{\infty}||U_s||_{L^{\infty}(\mathscr{A} \otimes \mathscr{A}^{op})}ds\right)^{1/2}
\end{equation*} where $L^{\infty}(\mathscr{A} \otimes \mathscr{A}^{op})$ 
is the norm defined by:  
\begin{equation*}
||\cdot||_{L^{\infty}(\mathscr{A} \otimes \mathscr{A}^{op})} = \lim_{p \rightarrow \infty} ||\cdot||_{L^p(\mathscr{A} \otimes \mathscr{A}^{op})}\end{equation*}
The completion of this space is denoted by $\mathscr{B}_{\infty}^a$. Recall also that, for fixed $t > 0$ and $U \in \mathscr{B}_{\infty}^a$, we have: 
\begin{equation*}
 \int_0^t U_s \sharp dX_s = \int_0^{\infty}U_s{\bf 1}_{[0,t]}(s)\sharp dX_s \end{equation*} where $X$ is a free additive Brownian motion and \begin{equation*}
 U_s \sharp dX_s := \sum_i A_s^i dX_s B_s^i\end{equation*}

Consider the process $(J_t := PY_tQY_t^{\star}P)_{t \geq 0}$, where $P, Q$ are two projections as in the definition above, $Y$ is a free multiplicative Brownian motion in  $(\mathscr{A}, \Phi)$. Recall that $Y$ satisfies the free SDE (cf \cite{Bia}): 
\begin{equation*} dY_t = i\,dX_t\,Y_t - \frac{1}{2} Y_t dt , \qquad Y_0 \in \mathscr{A}\end{equation*} where $(X_t)_{t \geq 0}$ is a free additive Brownian motion in  $(\mathscr{A}, \Phi)$.
By free It\^o's formula (\cite{Bia}, \cite{Bia1}, \cite{Kum}), we get :
\begin{align*}d(Y_tQY_t^{\star}) &= (dY_t)QY_t^{\star} + Y_tQ(dY_t^{\star}) + \Phi(Y_tQY_t^{\star}) dt \\&
= idX_t Y_tQY_t^{\star} -iY_tQY_t^{\star}dX_t - Y_tQY_t^{\star} +\theta dt\end{align*}
 since $X_t$ is self-adjoint. Thus, the free Jacobi process satisfies: \begin{align*}
dJ_t  &= P(dY_t\, Q\, Y_t^{\star})P\\
          &= iPdX_tY_tQY_t^{\star}P  - iPY_tQY_t^{\star}dX_tP - J_t\, dt  + \theta P dt \\
          &= iPY_t(Y_t^{\star}dX_tY_t)QY_t^{\star}P - iPY_tQ(Y_t^{\star}dX_tY_t)Y_t^{\star}P - J_t\, dt  + \theta P dt,\end{align*}
 since $\Phi$ is tracial and $Y_t$ is unitary (by definition). Next, let us recall the characterisation of the free Brownian motion (cf \cite{Bia2}, \cite{Cap}) which is the analogue of the classical L\'evy characterisation: Let $(\mathscr{A}_s,\, s \in [0,1])$ be an increasing family of von Neumann subalgebras in a non-commutative probability space $(\mathscr{A}, \Phi)$, and let $(Z_s = (Z_s^1, \hdots, Z_1^m); \, s \in [0,1])$ be an $m$-tuple of self-adjoint $(\mathscr{A}_s)$-adapted processes such that : \begin{enumerate}
 \item $Z$ is bounded and $Z_0 =0$.
 \item $\Phi(Z_t^i|\mathscr{A}_s) = Z_s^i$ for all $1 \leq i \leq m$.  
 \item $\Phi(|Z_t^i - Z_s^i|^4) \leq K(t-s)^2$ for some constant  $K$ and for all $1 \leq i \leq m$.
 \item For any $l, p \in \{1, \hdots, m\}$ and all $A,B \in \mathscr{A}_s$, one has: 
\begin{equation*}\Phi(A(Z_t^p -Z_s^p)B(Z_t^l - Z_s^l)) = {\bf 1}_{\{p=l\}}\Phi(A)\Phi(B)(t-s) + o(t-s),
\end{equation*} \end{enumerate}
then $Z$ is a $m$-dimensional free Brownian motion. Hence, we can see that :
\begin{lem}\label{Adnouna}The process $ (S_t) := (\int_0^tY_s^{\star}dX_sY_s)_{t \geq 0}$ is an $\mathscr{A}_t$ - free Brownian motion.\end{lem} 
{\it Proof} : Recall first that (cf \cite{Bia1}) for any two adapted biprocesses $N$ and $M$ belonging to $\mathscr{B}_{\infty}^{a}$: 
\begin{equation} \label{E1}\Phi\left\{\int_0^t N_s \sharp dX_s\, \int_0^t M_s \sharp dX_s \right\} = \int_0^t < N_s , M^{\star}_s > \,ds , \end{equation} 
where $X$ denotes a free Brownian motion in $(\mathscr{A},\Phi)$ and $< >$ is the inner product in $L^2(\mathscr{A} , \Phi) \otimes L^2(\mathscr{A} , \Phi)$ ( namely, if $N = a \otimes a'$ and $M = b \otimes b'$, then $M^{\star} = (b')^{\star} \otimes b^{\star}$ and $< N , M^{\star}> = \Phi(ab')\Phi(a'b)$). So , we have to verify the three conditions mentionned above. 
\\Note that for all $T > 0$, $Y_t^{\star} \otimes Y_t \mathbf{1}_{[0,T]} \in \mathscr{B}_{\infty}^{a}$ since $||Y_t|| = ||Y_t^{\star}|| = 1$ (Here, $|| \cdot ||$ denotes the algebra or $L^{\infty}(\mathscr{A}, \Phi)$ norm defined by $|| x || = \lim_{n \rightarrow \infty} \Phi(|x|^{2n})^{1/2n}$). Take $A, \, B \in \mathscr{A_s}$, then (using $(\ref{Adnouna})$ in the third line) : \begin{align*}
\Phi(A(S_t - S_s)B(S_t - S_s)) &= \Phi\left\{\int_s^t AY_r^{\star}dX_rY_r  \, \int_s^t BY_r^{\star}dX_rY_r \right\} \\
&= \Phi\left\{\int_s^t (AY_r^{\star}\otimes Y_r) \sharp dX_r  \, \int_s^t (BY_r^{\star}\otimes Y_r)\sharp dX_r \right\}\\
&= \int_s^t \Phi(Y_rAY_r^{\star}) \Phi(Y_rBY_r^{\star}) dr\\
&= \Phi(A)\Phi(B)(t - s), \end{align*} since $\Phi$ is tracial and $Y_t$ is unitary. Hence, the third condition is fullfilled. For the second one, we follow in the same way and use again the fact that $Y_t$ is unitary to get :  \begin{equation*} \Phi(|S_t - S_s|^4) \leq (t - s)^2 ,  \end{equation*}
The first one follows since $\left(\int_0^t Y_s^{\star}\, dX_s Y_s\right)_{t \geq 0}$ defines an $\mathscr{A}_t$ - martingale.$\hfill \square$\\ 
Thus, our free SDE is written: \begin{align*} dJ_t &=  iPY_t \,dS_t \,QY_t^{\star}P - iPY_tQ \,dS_t \,Y_t^{\star}P - J_t\, dt  + \theta P dt
\\& =   iPY_t (1 - Q)\,dS_t \,QY_t^{\star}P - iPY_tQ \,dS_t\,(1 - Q)Y_t^{\star}P - J_t\, dt  + \theta P dt\end{align*} where $1$ denotes the unit of $\mathscr{A}$. On the other hand, we can see that $P - J_t$ is positive. Indeed, since $Y_t$ is unitary, $Q^2 = Q$ and $P^2 = P$, one has : \begin{equation*}P - J_t = (PY_t - PY_tQ)(Y_t^{\star}P - QY_t^{\star}P) := C^{\star}C\end{equation*} Hence, using the following polar decompositions : \begin{eqnarray*}
QY_t^{\star}P &=& R_t \sqrt{J_t} \\
C = (1-Q)Y_t^{\star}P &=& V_t\sqrt{P - J_t}, \end{eqnarray*} 
our free SDE becomes : \begin{equation*}  dJ_t =   \sqrt{P - J_t} (iV_t^{\star}\,dS_t \,R_t) \sqrt{J_t} + \sqrt{J_t}((iR_t)^{\star} \,dS_t\, V_t)\, \sqrt{P - J_t} + \left(\theta P - J_t\right) dt  \end{equation*}  
\begin{nota}From the polar decomposition, we can see that : \begin{equation}\label{imp} \sqrt{J_t}R_t^{\star}V_t\sqrt{P-J_t} = PY_tQ(1-Q)Y_t^{\star}P = 0\end{equation}since $Q = Q^2$.\end{nota}
\begin{pro}\label{Bambola}
Suppose that $J_t$ and $P- J_t$ are injective operators in $P\mathscr{A}P$. Then, the following holds : 
\begin{itemize}
\item $R_tP =R_t \quad \textrm{and} \quad V_tP = V_t$.
\item $R_t^{\star}R_t = P  \quad \textrm{and} \quad V_t^{\star}V_t = P$. 
\item $PR_t^{\star} V_tP = PV_t^{\star}R_tP  = 0 $.
\end{itemize} \end{pro}  

{\it Proof}: Recall first that if $T$ is an operator in $\mathscr{A}$, then the support $E$ of $T$ is the orthogonal projection on $(ker \,T)^{\bot} = \overline{Im\, T^{\star}}$ and verifies $TE = T$ (cf Appendix III in \cite{Dix}). Furthermore, if we consider the polar decomposition of $T$, namely :\begin{equation*}
T =A\,|T| = A(T^{\star}T)^{1/2} ,\end{equation*}then, $E$ is also the support of $A$ and the latter is partially isometric, that is $A^{\star}A = E$ and  $A A^{\star} = F$  where $F$ is the support of $T^{\star}$. Thus, the two first assertions follow if we prove that $P$ is the support of both $J_t$ and $P-J_t$. Indeed, the injectivity of $J_t$ in $P\mathscr{A}P$ implies that $ker\,J_t =  ker\, P$. Thus, we claim that $P$ is the support of $J_t$ ($(ker \, P)^{\bot} = Im \, P$) and the same result holds for $P-J_t$. For the third, when $J_t$ and $P-J_t$ are invertible, it is obvious. Else, (\ref{imp}) is written in $P\mathscr{A}P$ : \begin{equation*}
0 = \sqrt{J_t}R_t^{\star}V_t\sqrt{P-J_t} = \sqrt{J_t}(PR_t^{\star}V_tP)\sqrt{P-J_t}\end{equation*}  Since both $J_t$ and $P-J_t$ are injective operators in $P\mathscr{A}P$, then :
\begin{equation*}
(PR_t^{\star}V_tP)\sqrt{P-J_t} = 0 \Rightarrow \sqrt{P-J_t}(PV_t^{\star}R_tP) = 0 \Rightarrow (PV_t^{\star}R_tP) = 0\end{equation*}  

\begin{cor} \label{Zizou}Under the same assumption of Proposition \ref{Bambola} , the process $(W_t)_{t \geq 0}$ defined by $W_t := (i/\sqrt{\Phi(P)})\int_0^t (PV_s^{\star} \otimes R_sP) \sharp dS_s$ is a $P\mathscr{A}_tP$ -complex free Brownian motion.\end{cor}
{\it Proof}: Let us first recall that a process $Z : \re_+ \rightarrow \mathscr{A}$ is a complex $(\mathscr{A}_t)$ -Brownian motion if it can be written $\displaystyle Z = \frac{X^1 + i X^2}{\sqrt 2}$, where $(X^1, X^2)$ is a $2$-dimensional $(\mathscr{A}_t)$ -free Brownian motion. Note also that $(iZ_t)_{t \geq 0}$ is still a complex $(\mathscr{A}_t)$ -Brownian motion since $(-X^2)$ is an  $(\mathscr{A}_t)$ -free Brownian motion. So, we shall show that :
\begin{equation*}\left(\tilde{W}_t := (1/\sqrt{\Phi(P)})\int_0^t (PV_s^{\star} \otimes R_sP)\sharp dS_s\right)_{t \geq 0}\end{equation*} is a $P\mathscr{A}_tP$ -complex free Brownian motion. To do this, it suffices to show that :
\begin{eqnarray*} 
X^1_t &=& \frac{\tilde{W}_t + \tilde{W}_t^{\star}}{\sqrt 2} = \frac{1}{\sqrt {2\Phi(P)}} \left(\int_0^t (PV_s^{\star} \otimes R_sP)\sharp dS_s + \int_0^t (PR_s^{\star} \otimes V_sP)\sharp dS_s\right) \\
X^2_t &=& \frac{\tilde{W}_t - \tilde{W}_t^{\star}}{\sqrt 2i} = \frac{1}{\sqrt {2\Phi(P)}i} \left(\int_0^t (PV_s^{\star} \otimes R_sP)\sharp dS_s - \int_0^t (PR_s^{\star} \otimes V_sP)\sharp dS_s\right) \end{eqnarray*} define two free $(\mathscr{A}_t)$ -free Brownian motions using the characterisation of free Brownian motion. We will do this for $X_1$. Note that, since $R_t$ and $V_t$ are partially isometric, then: 
\begin{equation*}(PV_t^{\star} \otimes R_tP\mathbf{1}_{[0, T]})_{t \geq 0} \quad \textrm{and} \quad (PR_t^{\star} \otimes V_tP\mathbf{1}_{[0, T]})_{t \geq 0} \in  B_{\infty}^a \, \forall \,T > 0
\end{equation*}  Hence, the first condition follows since $\left(\int_0^t (PV_s^{\star} \otimes R_sP) dS_s\right)_{t \geq 0}$ and $\left(\int_0^t (PR_s^{\star} \otimes V_sP) dS_s\right)_{t \geq 0}$ are $P\mathscr{A}_tP$-martingales. For $A, B \in \mathscr{A}_s$ and using (\ref{E1}), one has : 
\begin{align*}
&\tilde{\Phi}(PAP(X^1_t -X^1_s)PBP(X^1_t -X^1_s)) = \\& \frac{1}{2\Phi(P)}\tilde{\Phi}\left(\int_s^t (PAPV_u^{\star} \otimes R_uP)\sharp dS_u +  \int_s^t (PAPR_u^{\star} \otimes V_uP) \sharp dS_u\right)\\& \left(\int_s^t (PBPV_u^{\star} \otimes R_uP)\sharp dS_u +  \int_s^t (PBPR_u^{\star} \otimes V_uP)\sharp dS_u\right)
\\& = \frac{1}{2\Phi^2(P)}\Phi\left(\int_s^t (PAPV_u^{\star} \otimes R_uP)\sharp dS_u +  \int_s^t (PAPR_u^{\star} \otimes V_uP) \sharp dS_u\right)\\& \left(\int_s^t (PBPV_u^{\star} \otimes R_uP)\sharp dS_u +  \int_s^t (PBPR_u^{\star} \otimes V_uP)\sharp dS_u\right)
\\& = \frac{1}{2}\int_s^t [\tilde{\Phi} (PAPV_u^{\star}V_uP) \tilde{\Phi}(R_uPBPR_u^{\star}) +\tilde{\Phi} (PAPR_u^{\star}R_uP) \tilde{\Phi}(V_uPBPV_u^{\star})]du \\&
+ \frac{1}{2}\int_s^t [\tilde{\Phi} (PAPV_u^{\star}R_uP) \tilde{\Phi}(R_uPBPV_u^{\star}) +\tilde{\Phi} (PAPR_u^{\star}V_uP) \tilde{\Phi}(V_uPBPR_u^{\star})]du \\&
=\frac{1}{2}\int_s^t [\tilde{\Phi} (PAP) \tilde{\Phi}(PBP) +\tilde{\Phi} (PAP) \tilde{\Phi}(PBP)]du + \\&\frac{1}{2}\int_s^t [\tilde{\Phi} (APV_u^{\star}R_uP) \tilde{\Phi}(BPV_u^{\star}R_uP) +\tilde{\Phi} (APR_u^{\star}V_uP) \tilde{\Phi}(BPR_u^{\star}V_uP)]du  \\& = \tilde{\Phi}(PAP)\tilde{\Phi}(PBP)(t-s),\end{align*}  
since $PV_u^{\star}R_uP = PR_u^{\star}V_uP = 0$ and since $R_u$ and $V_u$ are partially isometric (Proposition \ref{Bambola}). Using the same arguments, we can see that the same result holds for $X^2$. Furthermore, one has :  \begin{align*}
&\tilde{\Phi}(PAP(X^1_t -X^1_s)PBP(X^2_t -X^2_s)) = \\& \frac{1}{2\Phi(P)}\tilde{\Phi}\left(\int_s^t (PAPV_u^{\star} \otimes R_uP)\sharp dS_u +  \int_s^t (PAPR_u^{\star} \otimes V_uP) \sharp dS_u\right)\\& \left(\int_s^t (PBPV_u^{\star} \otimes R_uP)\sharp dS_u -  \int_s^t (PBPR_u^{\star} \otimes V_uP)\sharp dS_u\right)\\& 
= \frac{1}{2}\int_s^t [-\tilde{\Phi} (PAPV_u^{\star}V_uP) \tilde{\Phi}(R_uPBPR_u^{\star}) +\tilde{\Phi} (PAPR_u^{\star}R_uP) \tilde{\Phi}(V_uPBPV_u^{\star})]du \\&
+ \frac{1}{2}\int_s^t [\tilde{\Phi} (PAPV_u^{\star}R_uP) \tilde{\Phi}(R_uPBPV_u^{\star}) - \tilde{\Phi} (PAPR_u^{\star}V_uP) \tilde{\Phi}(V_uPBPR_u^{\star})]du \\&
=\frac{1}{2}\int_s^t [-\tilde{\Phi} (PAP) \tilde{\Phi}(PBP) +\tilde{\Phi} (PAP) \tilde{\Phi}(PBP)]du + \\&\frac{1}{2}\int_s^t [\tilde{\Phi} (APV_u^{\star}R_uP) \tilde{\Phi}(BPV_u^{\star}R_uP) -\tilde{\Phi} (APR_u^{\star}V_uP) \tilde{\Phi}(BPR_u^{\star}V_uP)]du  = 0\end{align*}which finishes the proof. Substituting $R_t$ and $V_t$ by $R_tP$ and $V_tP$ and using Corollary 
(\ref{Zizou}), we get : 
\begin{teo}\label{th1}
Let $J_0$ such that $J_0$ and $P-J_0$ are injective operators in $P\mathscr{A}P$ and let $T := \inf\{s, \textrm{ker}(J_s) \neq ker P \quad \textrm{or} \quad \textrm{ker}(P - J_s) \neq ker P\}$. 
Then, for all $t < T$, we have
\begin{equation} 
\left\{\begin{array}{ll} dJ_t & =  \sqrt{\lambda \theta} \sqrt{P - J_t} dW_t \sqrt{J_t} + \sqrt{\lambda \theta}\sqrt{J_t}dW_t^{\star}\, \sqrt{P - J_t} + \left(\theta P - J_t\right) dt  \label{free}
 \\ J_0 & = PY_0QY_0^{\star}P
\end{array}\right.
\end{equation} where $(W_t)_{t \geq 0}$ is a $P\mathscr{A}P$-complex free Brownian motion.\end{teo} 
In the remainder of this paper, we will try to find the range of $(\lambda,\theta)$ ensuring the injectivity of both $J_t$ and $P-J_t$. We first investigate the stationary case then deal with the general setting.
\section{Free Jacobi process : the stationary case} 
In this section, we will give some interest in the particular case when $Y_0$ is Haar distributed, that is $\Phi(Y_0^k) = \delta_{k0}$. Then $Y_t$ remains Haar distributed for all $t > 0$. Thus, the law of $J_t$ does not depend on time and such a process is called a stationary free Jacobi process. We will compute the law of $J_t$ and  try to find values of $\lambda$ and $\theta$ such that there is no Dirac mass at $0$ and $1$, i. e, neither $0$, nor $1$ belongs to the discrete spectrum of $J_t$.  \\
Our method differs from the ones used by both Capitaine and Casalis (cf \cite{Cap}) and Collins ($P=Q$, cf \cite{Col}). It relies on Nica and Speicher's result concerning the compression by free projections.  Namely, authors considered $PaP$ for any  operator $a \in \mathscr{A}$ free with $P$ (cf \cite{Nica}, \cite{Spei}). This condition is fulfilled for 
$a=Y_tQY_t^{\star}$ since $Y_t$ is Haar unitary. In fact,  the following classical result holds (cf \cite{Hiai}):
\begin{lem}If $U$ is Haar unitary and $\mathscr{B}$ is a sub-algebra which is free with $U$, then, $\forall A, B \in \mathscr{B}$, $A$ and $UBU^{\star}$ are free.\end{lem}
Thus, denoting the law of $a$ in $(\mathscr{A},\Phi)$ by $\mu_a$, the law of $J_t$  in $(P\mathscr{A}P, \tilde{\Phi})$ is given by (cf \cite{Spei}): 
\begin{equation*}\mu = \boxplus^r \mu_{\alpha a} \end{equation*} where $\Phi(P) = \alpha$, $r = 1/\alpha$ and $\boxplus$ denotes the free additive convolution. In our case, $\alpha = \lambda \theta$ and $a = Y_tQY_t^{\star}$ so that for all $k \geq 1$, $\Phi(a^k) = \Phi(Q) = \theta$ since $\Phi$ is tracial and $Q$ is a projection. Thus, 
\begin{equation*} \mu_a = (1-\theta)\delta_0 + \theta \delta_1 \end{equation*} Next, recall that the  R-transform linearize the free additive convolution, so :
\begin{equation*} R(z) = r R_{\alpha a}(z) = R_a(\alpha z)\end{equation*} where the last equality follows from the definition of the R-transform and the multilinearity of free cumulants. Furthermore, \begin{equation*} G_a(z) = \frac{1}{z} + \sum_{k \geq 1}\frac{\Phi(a^k)}{z^{k+1}} = \frac{z+\theta - 1}{z(z-1)},\end{equation*}    
Hence, \begin{equation*} K_a(z) = \frac{z+1+ \sqrt{(z-1)^2+4\theta z}}{2z},\end{equation*} and finally, \begin{equation*}R_a(z) = K_a(z) - \frac{1}{z} = \frac{z-1+ \sqrt{(z-1)^2+4\theta z}}{2z}.\end{equation*}Consequently, we get : \begin{align*}R(z) &= \frac{\alpha z-1+ \sqrt{(\alpha z-1)^2+4\theta \alpha z}}{2\alpha z}
\\& = \frac{ z- r + \sqrt{(z-r)^2+4z/\lambda}}{2z},\end{align*} which implies that : 
\begin{equation*} K(z) = \frac{z +(2 - r) + \sqrt{(z-r)^2+4z/\lambda}}{2z}.\end{equation*}Thus, \begin{equation*}
G(z) = \frac{(2-r)z + (1/\lambda -1) + \sqrt{Az^2 - Bz + C}}{2z(z-1)},\end{equation*}where $A=r^2$, $B= 2(r + (r - 2)/\lambda)$ et $C = (1-1/\lambda)^2$.
Hence, the $J_t$ 's law is given by: \begin{equation*} \mu(dx) = a_0\delta_0(dx) + a_1\delta_1(dx) +  g(x)dx, \end{equation*} where \begin{eqnarray*}
a_0 &=& \lim_{y \rightarrow 0^+} -y\Im[G(iy)], \quad a_1 = \lim_{y \rightarrow 0^+} -y\Im[G(1 + iy)]\\
g(x) &=& \lim_{y \rightarrow 0^+} -\frac{1}{\pi}\Im[G(x + iy)] \qquad \textrm{for some} \, x \in (0,1),\end{eqnarray*} 
\begin{nota} The third equality holds whenever $\lim_{z \in D \rightarrow x}\Im(G(z)) = \Im(G(x))$ where $D$ is the upper half-plane and $x \in \re$ (cf \cite{Choi}). In fact, if we denote by \begin{equation*}
G(z) = \int \frac{dF(\zeta)}{z - \zeta} \end{equation*} the Cauchy transform of a distribution function $F$, then the following inversion formula holds (cf \cite{Hiai}, \cite{Spei}): \begin{equation*}
F(b) - F(a) = \lim_{y \rightarrow 0^+}-\frac{1}{\pi}\int_a^b \Im(G(x+iy)) dx \end{equation*} for any two continuity points $a, b$ of $F$ (weak convergence). Silverstein and Choi showed that if the limit above exists, then $F$ is differentiable and $dF$ has the density function with respect to the Lebesgue measure given by $F'(x) = lim_{y \rightarrow 0^+} -(1/\pi)  \Im(G(x+iy))$. See \cite{Choi} for more details.\end{nota}

From \cite{Capi}, we deduce that :  
\begin{equation*}\mu_{l,k}(dx) = \max(0, 1- l)\delta_0(dx) + \max(0, 1-k)\delta_1(dx) + g(x)\mathbf{1}_{[x_-,x_+]}dx \end{equation*}where
$l = 1/\lambda, k = (1-\theta /(\lambda \theta))$ and \begin{eqnarray*}
x_{\pm} & = &  \left(\sqrt{\theta(1 - \lambda \theta)} \pm  \sqrt{\lambda \theta(1 - \theta)}\right)^2 \\ 
g(x) & = & \frac{\sqrt{(x-x_-)(x_+ - x)}}{2\lambda \theta \pi x(1-x)} \end{eqnarray*}
Note that conditions $\lambda \in [0,1], 1/\theta \geq \lambda + 1$ implies that  $l \geq 1, k \geq 1$ so that $a_0 = b_0 =0$. 
Hence, for $\lambda = 1$, $x_- = 0$ and if $\theta = 1/\lambda +1$, then $x_+ = 1$. Consequently, we have :
 \begin{pro} $\forall \lambda \in ]0,1] , 1/\theta \geq \lambda +1\, (\theta \in ]0,1/2]$ for instance) , $J_t$ and $P-J_t$ are injective operators in the compressed space $P\mathscr{A}P$. For $\lambda \in ]0,1[ $ and $1/\theta > \lambda +1$, these operators are invertible in $P\mathscr{A}P$. Moreover, $J_t$ is a solution of (\ref{free}). \end{pro}

\begin{note}
1/In \cite{Capi}, authors omit the normalizing constant $\sqrt{A}/2\pi = (2\pi\lambda \theta)^{-1}$, however one can compute it as follows : since\begin{equation*}
\frac{\sqrt{(x_+ - x)(x- x_-)}}{x(1-x)} = \frac{\sqrt{(x_+ - x)(x- x_-)}}{x} + \frac{\sqrt{(x_+ - x)(x- x_-)}}{1-x}\end{equation*}
Then , the normalizing constant is given by  : \begin{align*}
 K^{-1} & = \int_{x_-}^{x_+}\frac{\sqrt{(x_+ - x)(x- x_-)}}{x}dx + \int_{1-x_+}^{1-x_-}\frac{\sqrt{(1- x_-  - x)(x- 1 + x_+)}}{x}dx  
\\& := I(x_-,x_+) + I(1-x_+, 1-x_-)\end{align*}
Moreover, using the variable change $u = (x-x_-)/(x_+ - x_-)$ and the integral representation of ${}_2F_1$ (cf \cite{Gra}), one has : \begin{align*}
I(x_-,x_+) &= \frac{(x_+ - x_-)^2}{x_+}B(\frac{3}{2},\frac{3}{2})\, {}_2F_1(1,\frac{3}{2},3,\frac{x_+-x_-}{x_+}) 
\\& \overset{(1)}{=} \frac{\pi}{4}\frac{(x_+ - x_-)^2}{x_+  + x_-}\, {}_2F_1(\frac{1}{2},1,2,\left(\frac{x_+-x_-}{x_+ + x_-}\right)^2) 
\\& \overset{(2)}{=} \frac{\pi}{2} (x_+ + x_- - 2\sqrt{x_+x_-}) \end{align*}
where in (1) we used (cf \cite{Gra}, p. 1070) : \begin{equation*} {}_2F_1(a,b,2b,z) = (1 - z/2)^{-a} \, {}_2F_1(a/2,(a+1)/2,b+1/2;(z/2-z)^2), \quad 0 < |z| < 1,\end{equation*}
and in (2), we used \begin{equation*} {}_2F_1(\frac{1}{2}, 1, 2; z) = 2\frac{1- \sqrt{1-z}}{z} , \qquad 0 < |z| < 1.\end{equation*} 
As a result,  $K^{-1}= \pi (1- \sqrt{x_+x_-} - \sqrt{(1-x_+)(1-x_-)}) = 2\pi\lambda \theta$. In the same way, setting $m_n(t) = \tilde{\Phi}(J_t^n), \forall n \geq 2$, one can see:  \begin{align*}
m_n -m_{n+1} &= K \int_{x_-}^{x_+}x^{n-1}\sqrt{(x_+ - x)(x- x_-)} \, dx  
\\& =\frac{K \pi(x_+ - x_-)^2}{8}x^{n-1}{}_2F_1\left(1-n,3/2,3; \frac{x_+ - x_-}{x_+},\right) \end{align*} 
\begin{equation*} m_1 - m_2 = \frac{K\pi}{8}(x_+-x_-)^2, \quad 1-m_1 = K  I(x_-,x_+).\end{equation*}
Thus, \begin{equation*}1- m_n = \frac{K\pi}{8}(x_+-x_-)^2\left[ 1+ \frac{1}{x_+}\left(1 + \sum_{\substack{k=0\\ k \neq 1}}^{n-1}x_+^{k}{}_2F_1\left(1-k,3/2,3; \frac{x_+ - x_-}{x_+},\right)\right)\right] \end{equation*}
One can also compute $m_n$ by expanding $(1-x)^{-1} = \sum_{k=0}^{\infty}x^k$ or use the integral representation of the Appell function $F_1$ (see ch. I in \cite{Ext}).\\ 

2/ One can compute $\tilde{\Phi}(\log(P-J_t))$ as follows : 
let $0 \leq z \leq 1$ and $\lambda \in ]0,1[, 1/\theta > \lambda +1$. Then, one can see that: \begin{equation*}
-\frac{d}{dz}\tilde{\Phi}(\log(P- zJ_t)) = -\frac{1}{z}\left(\frac{1}{z}G(\frac{1}{z}) - 1\right) = \frac{(1+ 1/\lambda)z -r + \sqrt{Cz^2 - Bz + A}}{2z(1-z)} \end{equation*}
Note that this derivative is well defined for $z=0$ and $z=1$.  It follows that : 
\begin{equation}\label{Haythem}  2\tilde{\Phi}(\log(P- J_t))  = -\int_0^1 \frac{(1+ 1/\lambda)z -r + \sqrt{Cz^2 - Bz + A}}{z(1-z)} dz\end{equation}
Note first that $Cz^2 - Bz+A > 0 \, \forall \lambda \in ]0,1[ , 1/\theta > \lambda + 1$ since $x_+,x_- \in ]0,1[$ are the roots of $Az^2 - Bz + C$ (so that $z < 1/x_+$) . In order to evaluate the integral in the right, we use the variable change $\sqrt{A}(1 - uz) = \sqrt{Cz^2 - Bz + A}$, which gives : \begin{equation*}
z = \frac{2Au - B}{Au^2- C}, \, \, 1-z = \frac{Au^2 - 2Au + B-C}{Au^2-C}, \, \, dz = -2A\frac{Au^2 - Bu + C}{(Au^2 - C)^2}du. \end{equation*}
Moreover, since $A - B + C = A(1- \theta(\lambda+1))^2 \geq 0$ and $\theta(1 + \lambda) \leq 1$, then the roots of $Au^2 - 2Au + B - C = 0$ are given by :
$u_{\pm} = 1\pm (1-\theta(\lambda + 1))$. On the other hand, $B/2A = (1/2)(x_+ + x_-) = \theta(\lambda+1-2\lambda \theta)$. Hence our expression factorizes to : 
\begin{align*}
\tilde{\Phi}(\log(P- J_t))   & = -\frac{1}{\sqrt A} \int_{B/2A}^{\theta(\lambda+1)} \frac{(u - \theta(\lambda+1))(Au^2 - Bu +C)}{(u^2 - C/A)(Au^2 - 2Au + B - C)} du  
\\& = -\frac{1}{\sqrt A} \int_{B/2A}^{\theta(\lambda+1)} \frac{(Au^2 - Bu +C)}{(u^2 - C/A)(u - u_+)} du  
\\& =   \int_{B/2A}^{\theta(\lambda+1)} \frac{C_1}{u - \theta(1-\lambda)} + \frac{C_2}{u + \theta(1-\lambda)}   + \frac{C_3}{u - u_+}  du \end{align*}
for some constants $C_1, C_2, C_3$ depending on both $\lambda, \theta$, given by : \begin{eqnarray*}
C_1 & = & 1 , \qquad C_2 = 1/\lambda, \qquad   C_3  = \frac{1-\theta(\lambda+1)}{\lambda \theta}\end{eqnarray*}  
Thus, one gets : \begin{align*}
\tilde{\Phi}(\log(P- J_t)) &= -\left[C_1\log(u-\theta(1-\lambda)) + C_2\log(u+\theta(1-\lambda)) + C_3\log(u_+ - u)\right]_{B/2A}^{\theta(\lambda +1)}
\\& = \log(1-\theta) + \frac{1}{\lambda} \log(1-\lambda \theta) - C_3 \log\left[\frac{(1-\theta(\lambda+1))}{1 - \theta(\lambda + 1) + \lambda \theta^2}\right]
\\& = (1+C_3)\log(1-\theta) + (\frac{1}{\lambda}+C_3) \log(1-\lambda \theta) - C_3 \log(1-\theta(\lambda+1))
\\& = \frac{(1-\theta)\log(1-\theta)  + (1-\lambda \theta)\log(1-\lambda\theta) - (1-\theta(\lambda+1))\log(1-\theta(\lambda+1))}{\lambda \theta} \end{align*}
Note that the result extends for all $\lambda \in ]0,1], \, 1/\theta \geq \lambda + 1$. Since $P-J$ is still a $FJP(\lambda \theta /(1-\theta),1-\theta)$, then : \begin{equation*}
\tilde{\Phi}(\log(J_t)) = \frac{\theta\log\theta  + (1-\lambda \theta)\log(1-\lambda\theta) - \theta(1-\lambda)\log(\theta(1-\lambda))}{\lambda \theta} \end{equation*}
which agrees with Rouault's result (cf \cite{Rouault}). \\
3/ In \cite{Yan}, Doumerc finds that for $p(m) \geq m+1$ and $q(m) \geq m+1$ where $q(m) = d(m) - p(m)$, the Jacobi process satisfies the following SDE: \begin{equation*}
dJ_t =  \sqrt{I_m - J_t} dB_t \sqrt{J_t} + \sqrt{J_t}dB_t^T\, \sqrt{I_m - J_t} + \left(p(m)I_m - (p(m) + q(m)) J_t\right) dt\end{equation*} where $(B_t)$ is a real $m \times m$ Brownian matrix.
The complex version of this process verifies:  \begin{equation*}
dJ_t =  \sqrt{I_m - J_t} dB_t \sqrt{J_t} + \sqrt{J_t}dB_t^{\star}\, \sqrt{I_m - J_t} + \left(p(m)I_m - (p(m) + q(m)) J_t\right) dt\end{equation*} where $(B_t)$ is a complex $m \times m$ Brownian matrix. Heuristically, if we consider the ratio $dJ_t/ (d(m))$ and let $m$ go to infinity, then this SDE converges weakly (up to a constant) to its free counterpart, since normalized complex Brownian matrix converges in distribution to the free complex Brownian motion. It is also worthnoting that conditions  $p(m) \geq m+1$ and $q(m) \geq m+1$ are in agreement with $\lambda \in [0,1[$ and $1/\theta \geq \lambda +1$. \\
4/ From a combinatorial point of view, the result of Nica and Speicher reads :
\begin{equation*}
\Phi(J_t^n) = \sum_{\pi \in NC(n)} k_{\pi}(P,\dots,P)\theta^{n+1-|\pi|}\end{equation*} 
where $NC(n)$ denotes the set of non-crossing partitions of $\{1, \dots,n\}$, $|\pi|$ is the cardinality of $\pi$ and $k_{\pi}$ is the corresponding mixed cumulant (see \cite{Spei} for more details). \end{note}

\section{Free Jacobi Process : the general setting.}
In this section, we will suppose that $\lambda \leq 1$ and $1/\theta \geq \lambda +1$. Let  $J_0$ in $P\mathscr{A}P$ such that $0 < J_0 < P$, means $J_0$ and $P-J_0$ are invertible in $P\mathscr{A}P$. Then, by continuity of paths, the result of theorem \ref{th1} holds for $t < T$, that is : 
\begin{equation*}\left\{\begin{array}{ll}
dJ_t = &  U_t \sharp dX_t  + V_t  \sharp dY_t + \left(\theta P - J_t\right) dt \\ J_0 = & PY_0QY_0^{\star}P ,\qquad 0< J_0 < P\end{array}\right. \end{equation*}
where \begin{eqnarray*}  W_t & = & \frac{X_t + iY_t}{\sqrt{2}}  \\
U_t & = & \sqrt{\frac{\lambda \theta}{2}}(\sqrt{P-J_t} \otimes \sqrt{J_t} + \sqrt{J_t} \otimes \sqrt{P - J_t}) = \sum_{i=1}^2A_t^i \otimes B_t^i \\ 
V_t & = & i\sqrt{\frac{\lambda \theta}{2}}(\sqrt{P-J_t} \otimes \sqrt{J_t} - \sqrt{J_t} \otimes \sqrt{P - J_t}) = \sum_{i=1}^2C_t^i \otimes D_t^i\end{eqnarray*} and $X$ and $Y$ are two 
free $P\mathscr{A}_tP$-free-Brownian motions. Now, let us recall that for any operator $Z$, we set ({cf \cite{Bia1}): 
\begin{eqnarray*} \partial Z^n & = & \sum_{k = 0}^{n-1}Z^k \otimes Z^{n-k-1} \\
\Delta_U(Z^n) & = & 2\sum_{i,j}\sum_{\substack{k,l \geq 0 \\ k+l \leq n-2}} Z^kA_t^iB_t^j Z^{n-k-l-2}\tilde{\Phi}(B_t^iZ^lA_t^j) \end{eqnarray*}where $U= \sum_i A^i \otimes B^i$ is an adapted bi-process. 
\begin{pro}
\label{S}
Let $X, Y$ be two free free-Brownian motions, $U, V$ be two adapted integrable bi-processes and $K$ an adapted process in $P\mathscr{A}P$. Let 
\begin{equation*}
dM_t = U_t\sharp dX_t  + V_t \sharp dY_t + K_tdt \end{equation*}then, for every polynomial $R$, we have : \begin{align*}
dR(M_t) &= (\partial R(M_t) \sharp U_t) \sharp dX_t + (\partial R(M_t) \sharp V_t) \sharp dY_t +  \partial R(M_t) \sharp K_t dt
\\& +\frac{1}{2}(\Delta_{U}R(M_t) + \Delta_{V}R(M_t))dt  \end{align*} \end{pro}

{\it Proof} : When $V = 0\otimes 0$, this is the free It\^o's formula stated in \cite{Bia1}. By linearity, it suffices to prove the formula for monomials. To do this, we shall proceed by induction. Hence, we assume that: \begin{align*}
dM_t^n & = \partial M_t^n \sharp (U_t \sharp dX_t) + \partial M_t^n \sharp (V_t \sharp dY_t) + \partial M_t^n \sharp K_t dt + \\& \frac{1}{2}(\Delta_{U}M_t^n + \Delta_{V}M_t^n)dt 
\end{align*}By free integration by parts formula (\cite{Bia}), we have:
\begin{align*}dM_t^{n+1} &= d(M_tM_t^n) = dM_t \, M_t^n + M_t \, dM_t^n + (dM_t)(dM_t^n)
\\& = \left((1 \otimes M_t^n + M_t\partial M_t^n) \sharp U_t\right) \sharp dX_t + \left((1 \otimes M_t^n + M_t\partial M_t^n) \sharp V_t\right) \sharp dY_t \\& 
+(1 \otimes M_t +  M_t\partial M_t^n) \sharp K_t dt    + \frac{1}{2}M_t(\Delta_{U}M_t^n + \Delta_{V}M_t^n)dt + (dM_t)(dM_t^n)\end{align*} 

On the other hand, we can easily see that : \begin{equation*}1 \otimes M_t^n + M_t\partial M_t^n = 1 \otimes M_t^n + \sum_{k=1}^{n}M_t^k \otimes M_t^{n-k} = \partial M_t^{n+1},\end{equation*}Then, using the fact that $(dX)(dY) =0$ by the freeness of $X$ and $Y$, we get : \begin{align*}
(dM_t)(dM_t^n) =  \sum_{i,j}\sum_{l=0}^{n-1}A_t^iB_t^j M_t^{n-l-1} \tilde{\Phi}(B_t^iM_t^lA_t^j) +  \sum_{i,j}\sum_{l=0}^{n-1}C_t^iD_t^j M_t^{n-l-1} \tilde{\Phi}(D_t^iM_t^lC_t^j) \end{align*}
Moreover, \begin{align*} M_t\Delta_U(M_t^n) & =  2\sum_{i,j}\sum_{\substack{k,l \geq 0\\  k+l \leq n-2}} M_t^{k+1} A_t^iB_t^j M_t^{n-k-l-2}\tilde{\Phi}(B_t^iM_t^lA_t^j) \\&
= 2\sum_{i,j}\sum_{l=0}^{n-2}\sum_{k=1}^{n-l-1} M_t^{k} A_t^iB_t^j M_t^{n-k-l-1}\tilde{\Phi}(B_t^iM_t^l A_t^j)\end{align*}and the same holds for \begin{equation*} 
M_t\Delta_V(M_t^n) = 2\sum_{i,j}\sum_{l=0}^{n-2}\sum_{k=1}^{n-l-1} M_t^{k} C_t^iD_t^j M_t^{n-k-l-1}\tilde{\Phi}(C_t^iM_t^lD_t^j)\end{equation*}Consequently, we get : 
\begin{equation*}\frac{1}{2}M_t(\Delta_V(M_t^n)  + \Delta_U(M_t^n)) +  (dM_t)(dM_t^n) = \frac{1}{2}(\Delta_V(M_t^{n+1})  + \Delta_U(M_t^{n+1}))\end{equation*}which completes the proof. $\hfill \blacksquare$ 

\subsection{A recurrence formula for free Jacobi moments}
\begin{cor}\label{mom}
Let $m_n(t) := \tilde{\Phi}(J_t^n)$ for $n \geq 2$ and $t < T$. Then, we have the following recurrence relation: \begin{equation*}
m_n(t) = m_n(0) - n\int_0^tm_n(s) ds + n\theta\int_0^tm_{n-1}(s) ds + \lambda \theta n\sum_{k=0}^{n-2}\int_0^tm_{n-k-1}(s)(m_k(s) - m_{k+1}(s))ds
 \end{equation*}
or equivalently, \begin{equation*}\frac{dm_n(t)}{dt} =  - nm_n(t)  + n\theta m_{n-1}(t)  + \lambda \theta n\sum_{k=0}^{n-2}m_{n-k-1}(t)(m_k(t) - m_{k+1}(t)) \end{equation*}\end{cor}
{\it Proof}: Using Proposition (\ref{S}), we get: \begin{equation*} dJ_t^n = \textrm{martingale} + \sum_{k=0}^{n-1}J_t^k(\theta P - J_t)J_t^{n-k-1}dt  + \frac{1}{2}(\Delta_{U}(J_t^n) + 
\Delta_{V}(J_t^n))dt \end{equation*} Next, we compute \begin{align*}
\Delta_U(J_t^n) & = 2\sum_{i,j = 1}^2\sum_{\substack{k,l \geq 0 \\k+l \leq n-2}}J_t^kA_t^iB_t^jJ_t^{n-k-l-2}\tilde{\Phi}(B_t^iJ_t^lA_t^j) \\&
\overset{(1)}{=}2\sum_{i,j = 1}^2\sum_{\substack{k,l \geq 0 \\k+l \leq n-2}}A_t^iB_t^jJ_t^{n-l-2}\tilde{\Phi}(B_t^iJ_t^lA_t^j)\\&
=2\sum_{i,j = 1}^2\sum_{l = 0}^{n-2}(n-l-1)A_t^iB_t^jJ_t^{n-l-2}\tilde{\Phi}(B_t^iJ_t^lA_t^j)\\&
=2\sum_{i,j = 1}^2\sum_{l = 0}^{n-2}(l+1)A_t^iB_t^jJ_t^l\tilde{\Phi}(B_t^iJ_t^{n-l-2}A_t^j)\\&
\overset{(2)}{=}2\sum_{l=0}^{n-2}(l+1)[\frac{\lambda \theta}{2}(P-J_t)J_t^l \tilde{\Phi}(J_t^{n-l-1}) + \frac{\lambda\theta}{2} J_t^{l+1}\tilde{\Phi}(J_t^{n-l-2}(P-J_t)) \\&
+\lambda \theta \sqrt{P-J_t}\sqrt{J_t}J_t^l\tilde{\Phi}(J_t^{n-l-2}\sqrt{J_t}\sqrt{P-J_t})]  
\end{align*}where in both $(1)$ and $(2)$, we used the fact that $A_t^i, B_t^j$ and $J_t$ commute $\forall \, i, j \in \{1,2\}$. Similarly, we get : \begin{align*}
\Delta_V(J_t^n) &= 2\sum_{l=0}^{n-2}(l+1)[\frac{\lambda \theta}{2}(P-J_t)J_t^l \tilde{\Phi}(J_t^{n-l-1}) + \frac{\lambda\theta}{2} J_t^{l+1}\tilde{\Phi}(J_t^{n-l-2}(P-J_t)) \\&
- \lambda \theta \sqrt{P-J_t}\sqrt{J_t}J_t^l\tilde{\Phi}(J_t^{n-l-2}\sqrt{J_t}\sqrt{P-J_t})]\end{align*} Thus, we have: \begin{equation*}
\frac{1}{2}(\Delta_U(J_t^n) + \Delta_V(J_t^n)) = \lambda \theta \sum_{l=0}^{n-2} (l+1) [(P-J_t)J_t^l \tilde{\Phi}(J_t^{n-l-1}) + J_t^{l+1}\tilde{\Phi}(J_t^{n-l-2}(P-J_t))]\end{equation*}
Taking the expectation, it yields 
\begin{align*} &\frac{1}{2}\tilde{\Phi}(\Delta_U(J_t^n) + \Delta_V(J_t^n)) = \lambda \theta \times\\&
\left(\sum_{l=0}^{n-2} (l+1) [\tilde{\Phi}((P-J_t)J_t^l) \tilde{\Phi}(J_t^{n-l-1})] + \sum_{l=0}^{n-2}(l+1)[\tilde{\Phi}(J_t^{l+1})\tilde{\Phi}(J_t^{n-l-2}(P-J_t))]\right)\\& 
=  \lambda \theta \times\left(\sum_{l=0}^{n-2} (l+1) [\tilde{\Phi}((P-J_t)J_t^l) \tilde{\Phi}(J_t^{n-l-1})] + \sum_{l=0}^{n-2}(n-l-1)[\tilde{\Phi}(J_t^{n-l-1})\tilde{\Phi}(J_t^{l}(P-J_t))]\right)
\\& = n\lambda \theta \sum_{l=0}^{n-2} [\tilde{\Phi}((P-J_t)J_t^l) \tilde{\Phi}(J_t^{n-l-1})]  \\&
 = n\lambda \theta \sum_{l=0}^{n-2}[m_{n-l-1}(t)(m_l(t) - m_{l+1}(t)]\end{align*}
 Furthermore, since $P-J_t$ and $J_t$ commute, we can easily see that : \begin{align*}
\tilde{\Phi}\left(\sum_{k=0}^{n-1}J_t^k(\theta P - J_t)J_t^{n-k-1}\right) = n\tilde{\Phi}(J_t^{n-1}(\theta P-J_t)) \end{align*} and the proof is complete. $\hfill \square$
 
 \begin{note} 1/ In the last proof, we interwined integrals and expectation since powers of $J_t$ are still positive self-adjoint operators.\\
 2/ For $n=1$, we obviously have: \begin{align*}\tilde{\Phi}(J_t) = \tilde{\Phi}(J_0) + \theta t  - \int_0^t \tilde{\Phi}(J_s)ds \end{align*} which implies that : 
$m_1(t) = (\tilde{\Phi}(J_0)-\theta)e^{-t} + \theta$ .  We can also deduce this result from the definition of the free Jacobi process, using moments of $Y_t$ (cf \cite{Bia}) and the following formula : 
 \begin{equation*} \tau(a_1b_1a_2b_2) = \tau(a_1a_2)\tau(b_1)\tau(b_2) +\tau(b_1b_2)\tau(a_1)\tau(a_2) - \tau(a_1)\tau(a_2)\tau(b_1)\tau(b_2),\end{equation*}where $\{a_1, a_2\}$ and $\{b_1, b_2\}$ are free families in some non commutative probability space $(\mathscr{B}, \tau)$. \\
 3/ In the stationary case and for $\lambda =1$, this formula reads : \begin{equation*} 
 (1-\theta)m_n = \theta \sum_{k=0}^{n-1}m_{n-k-1}(m_k - m_{k+1}) \end{equation*}
 For  $\theta = 1/2$, one has : \begin{equation*}
m_n = \sum_{k=0}^{n-1}m_{n-k-1}(m_k - m_{k+1}) \end{equation*}since $m_0 = 1$. Note also that $J_t$ is a Beta random variable $B(1/2,1/2)$, it follows that:  
\begin{align*} 
\frac{\Gamma(n+1/2)}{n!} = m_n = \frac{1}{2\sqrt{\pi}}\sum_{k=0}^{n-1}\frac{\Gamma(n-k-1/2)\Gamma(k+1/2)}{(n-k-1)!(k+1)!}
 = \sum_{k=1}^{n}\frac{\Gamma(n-k+1/2)\Gamma(k-1/2)}{2\sqrt{\pi}(n-k)! k!} \end{align*}
Using duplication formula, we get : \begin{align*} \frac{(2n)!}{(n!)^2} & = 2 \sum_{k=1}^{n}\frac{(2n-2k)!(2k-2)!}{((n-k)!)^2k!(k-1)!}
\\& = \sum_{k=1}^{n}\frac{(2n-2k)!(2k-2)!}{((n-k)!)^2k!(k-1)!} + \sum_{k=1}^{n}\frac{(2n-2k)!(2k-2)!}{(n-k)!(n-k+1)!((k-1)!)^2} 
\\& = (n+1) \sum_{k=1}^{n}\frac{(2n-2k)}{(n-k+1)((n-k)!)^2} \frac{(2k-2)!}{k((k-1)!)^2}\end{align*}  
Recall that the Catalan numbers are given by : $D_p = 1/(p+1)\binom{2p}{p}$. 
Thus, our recurrence formula  is equivalent to $D_n = \sum_{k=1}^n D_{n-k}D_{k-1}$. \end{note}
\begin{pro} \label{inj}If $J_0$ is invertible, then for all $\lambda \in ]0,1]$ , $1/\theta \geq 1+ \lambda$ and $t \geq 0$,  $P-J_t$ and $J_t$ are injective operators .\end{pro}
{\it Proof}: It is known that for a self - adjoint operator $a \in P\mathscr{A}P$ such that $0 < ||a|| < P$, we have: \begin{equation*}
\log(P-a) = -\sum_{n=1}^{\infty} \frac{a^n}{n}\end{equation*}
Since $\tilde{\Phi}(a^n) = \int x^n\mu(dx)$ for a positive compact supported measure $\mu$, we get :
\begin{equation*}\tilde{\Phi}(\log(P-a)) = -\sum_{n=1}^{\infty}\frac{\tilde{\Phi}(a^n)}{n} = -\tilde{\Phi}(a) - \sum_{n=2}^{\infty}\frac{\tilde{\Phi}(a^n)}{n}\end{equation*}Thus, substituting moments of $J_t$, one has (cf Corollary (\ref{mom})):   

\begin{align*} \tilde{\Phi}&(\log(P-J_t)) = \tilde{\Phi}(\log(P-J_0))  - \theta t + \int_0^t \tilde{\Phi}\left(\sum_{n=1}^{\infty}J_s^n\right)ds - \theta  \int_0^t \tilde{\Phi}\left(\sum_{n=1}^{\infty}J_s^{n}\right)ds \\& - \lambda \theta \int_0^t \sum_{n=2}^{\infty}\sum_{k=0}^{n-2}\tilde{\Phi}(J_s^{n-1-k})\tilde{\Phi}(J_s^k(P-J_s))ds \\&
= \tilde{\Phi}(\log(P-J_0))  - \theta t + (1-\theta)\int_0^t\tilde{\Phi}(J_s(P-J_s)^{-1})ds -  \lambda \theta \int_0^t\sum_{k=0}^{\infty}\sum_{n=1}^{\infty}\tilde{\Phi}(J_s^n)\tilde{\Phi}(J_s^k(P-J_s))ds \\&
= \tilde{\Phi}(\log(P-J_0))  - \theta t + (1-\theta)\int_0^t\tilde{\Phi}(J_s(P-J_s)^{-1})ds -  \lambda \theta\int_0^t\tilde{\Phi}(J_s(P-J_s)^{-1})\tilde{\Phi}(P)ds\\&
= \tilde{\Phi}(\log(P-J_0))  - \theta t + (1-\theta - \lambda \theta)\int_0^t\tilde{\Phi}(J_s(P-J_s)^{-1})ds \\&
= \tilde{\Phi}(\log(P-J_0))  - (1-\lambda \theta)t + (1-\theta - \lambda \theta)\int_0^t\tilde{\Phi}((P-J_s)^{-1})ds \end{align*}
Hence, we deduce that, if $\lambda \in ]0,1]$ and $1/\theta \geq 1+\lambda$ -which ensures that $1 - \theta - \lambda \theta \geq 0$ and $1-\lambda \theta \geq 0$- and if $J_0$ is an invertible operator, then: \begin{equation*}
 \tilde{\Phi}(\log(P-J_t)) + (1- \lambda \theta)t \geq \tilde{\Phi}(\log(P-J_0)) > -\infty  \qquad \forall \, t < T \end{equation*} which gives the injectivity of $P-J_t, \, \forall t \geq 0$. The second assertion follows since one can easily see from (\ref{free}) that , if $J$ is a free Jacobi process with parameters $FJP(\lambda, \theta)$ starting at $0 < J_0 < P$, then $P-J$ is still a free Jacobi process with parameters $(\lambda \theta /1-\theta, 1-\theta)$ starting at $0 < P-J_0 < P$, and since $\frac{1}{\theta} \geq \lambda + 1  \Rightarrow  \frac{\lambda \theta}{1 - \theta} \leq 1$ and $\lambda \leq 1 \Rightarrow (\lambda -1)\theta \leq 0 \Rightarrow  (\lambda \theta)/(1-\theta) + 1 = (\theta(\lambda-1) + 1)/(1-\theta) \leq 1/(1-\theta)$.

\begin{cor} Under the same conditions of proposition (\ref{inj}), the $FJP(\lambda, \theta)$  satisfies for all $t \geq 0$ the following SDE: \begin{equation*}dJ_t =  \sqrt{\lambda \theta} \sqrt{P - J_t} dW_t \sqrt{J_t} + \sqrt{\lambda \theta}\sqrt{J_t}dW_t^{\star}\, \sqrt{P - J_t} + \left(\theta P - J_t\right) dt\end{equation*} 
 where $W$ is a complex free Brownian motion.\end{cor}

\subsection{Free martingales polynomials} 
In this paragraph, we consider a stationary $FJP(1,1/2)$ starting at $J_0$, the law of which is the Beta law $B(1/2,1/2)$.
Recall that that a $\mathscr{A}_t$-adapted free process $(X_t)_{t \geq 0}$ is a $\mathscr{A_t}$-free martingale if and only if $\Phi(X_t|\mathscr{A}_s) = X_s$ (cf \cite{Ans},\cite{Bia}, 
\cite{Bia3}). 
\begin{pro}
Let $\mathscr{J}_t$ denotes the von Neumann subalgebra generated by $(J_s, \, s \leq t)$ and let $0 < r < 1 $. Then,  the process $R_t := ((1+re^t)P- 2re^t J_t)((1+re^t)^2 P - 4re^tJ_t)^{-1})_{t < -\ln r}$ is a $\mathscr{J}_t$-free martingale. 
 \end{pro}
{\it Proof :}  $R_t$ can be written as : 
\begin{align*} 
R_t &= \left[\frac{P}{1+ re^t} - 2\frac{re^t}{(1+re^t)^2}J_t\right] \left[P - \frac{4re^t}{(1+re^t)^2}J_t\right]^{-1}
\\&  =  \left[\frac{1-re^t}{2(1+ re^t)}P + \frac{1}{2}\left(P - \frac{4re^t}{(1+re^t)^2}J_t\right)\right]\left[P - \frac{4re^t}{(1+e^t)^2}J_t\right]^{-1} 
:=  \frac{1-re^t}{2(1+ re^t)}H_t + \frac{P}{2}
\end{align*} where
\begin{equation*}
H_t = \left[P - \frac{4re^t}{(1+re^t)^2}J_t\right]^{-1} = \sum_{n \geq 0}\frac{(4re^t)^n}{(1+re^t)^{2n}}J_t^n
\end{equation*} since $4re^t < (1+re^t)^2$ and $0 \leq J_t \leq P$ for all $t > 0$. It follows that : 
\begin{equation*}
2dR_t  =  \frac{1-re^t}{1+ re^t}dH_t - \frac{2re^t}{(1+re^t)^2}H_t dt 
\end{equation*}
On the other hand, one has for $1 \leq l \leq n-1$ and $n  \geq 2$ : 
\begin{equation*}
\tilde{\Phi}(J_t^{n-l}) = \frac{\Gamma(n-l+1/2)}{\sqrt{\pi}(n-l)!}, \qquad  \tilde{\Phi}(J_t^{n-l-1}(P-J_t)) = \frac{\Gamma(n-l-1/2)}{2\sqrt{\pi}(n-l)!}
\end{equation*}

From the proof of Proposition (\ref{S}), we deduce that, for all $t > 0$ :     
\begin{align*} 
&dJ_t^n = M_t + n(\frac{P}{2} - J_t)J_t^{n-1}dt  + \frac{1}{2}(\Delta_{U}(J_t^n) + \Delta_{V}(J_t^n))
\\&=  M_t + \left[n(\frac{P}{2} - J_t)J_t^{n-1}+ \sum_{l=1}^{n-1}l[(P-J_t)J_t^{l-1} \tilde{\Phi}(J_t^{n-l}) + J_t^l  \tilde{\Phi}(J_t^{n-l-1}(P-J_t))] \right]dt
\\& =  M_t + \left[n(\frac{P}{2} - J_t)J_t^{n-1} + \frac{1}{2\sqrt{\pi}} \sum_{l=1}^{n-1} l\left[(P-J_t)J_t^{l-1}\frac{\Gamma(n-l+1/2)}{(n-l)!} 
+ \frac{J_t^l }{2}\frac{\Gamma(n-l-1/2)}{(n-l)!}\right] \right]dt 
\\& = M_t+ \left[n(\frac{P}{2} - J_t)J_t^{n-1}+ \frac{1}{2\sqrt{\pi}} \left[\sum_{l=1}^{n-1} l \frac{\Gamma(n-l+1/2)}{(n-l)!}J_t^{l-1}  - \sum_{l=1}^{n-1} l(n-l-1) \frac{\Gamma(n-l-1/2)}{(n-l)!} J_t^l \right]\right]dt
\\& =  M_t +\left[ n(\frac{P}{2} - J_t)J_t^{n-1} + \frac{1}{2\sqrt{\pi}}\left[\sum_{l=1}^{n-1} l \frac{\Gamma(n-l+1/2)}{(n-l)!}J_t^{l-1}  
- \sum_{l=2}^{n-1} (l-1)(n-l) \frac{\Gamma(n-l+1/2)}{(n-l+1)!} J_t^{l-1} \right]\right]dt
\\& =  M_t + \left[n(\frac{P}{2} - J_t)J_t^{n-1} + \frac{1}{2\sqrt{\pi}}\left[\sum_{l=2}^{n-1} n \frac{\Gamma(n-l+1/2)}{(n-l+1)!}J_t^{l-1} + \frac{\Gamma(n-1/2)}{(n-1)!} P \right] \right]dt
\\& =  M_t + \left[n(\frac{P}{2} - J_t)J_t^{n-1} + \frac{n}{2\sqrt{\pi}}\sum_{l=1}^{n-1}  \frac{\Gamma(n-l+1/2)}{(n-l+1)!}J_t^{l-1}\right]dt
\\& = M_t - nJ_t^n dt +  \frac{n}{2\sqrt{\pi}}\sum_{l=1}^{n}  \frac{\Gamma(n-l+1/2)}{(n-l+1)!}J_t^{l-1} \,dt \end{align*} 
where $M_t$ stands for the martingale part. Note that this holds for $n=1$. Thus : 

\begin{align*}
&FV(dH_t) = \sum_{n \geq 0}\frac{(4re^t)^n}{(1+re^t)^{2n}}FV(dJ_t^n) + \frac{1-re^t}{1+re^t}\sum_{n \geq 0}\frac{n(4re^t)^n}{(1+re^t)^{2n}}J_t^n dt
\\& = -\frac{2re^t}{1+re^t}\sum_{n \geq 0}\frac{n(4re^t)^n}{(1+re^t)^{2n}}J_t^n dt + \frac{1}{2\sqrt{\pi}}\sum_{n \geq 0}\frac{n(4re^t)^n}{(1+re^t)^{2n}} \sum_{l=1}^{n}  \frac{\Gamma(n-l+1/2)}{(n-l+1)!}J_t^{l-1} \,dt 
\\& = -\frac{2re^t}{1+re^t}\sum_{n \geq 0}\frac{n(4re^t)^n}{(1+re^t)^{2n}}J_t^n dt + \frac{1}{2}\sum_{l \geq 1}\frac{(4re^t)^l}{(1+re^t)^{2l}} 
\sum_{n \geq 0}  \frac{(n+l)(1/2)_n}{(n+1)!}\frac{(4re^t)^n}{(1+re^t)^{2n}}J_t^{l-1} \,dt 
\\& = -\frac{2re^t}{1+re^t}\sum_{n \geq 0}\frac{n(4re^t)^n}{(1+re^t)^{2n}}J_t^n dt + \frac{2re^t}{(1+re^t)^2}\sum_{l \geq 0}\frac{(4re^t)^l}{(1+re^t)^{2l}} \sum_{n \geq 0}  \frac{(n+l+1)(1/2)_n}{(n+1)!}\frac{(4re^t)^n}{(1+re^t)^{2n}}J_t^{l} \,dt 
\\& =  -\frac{2re^t}{1+re^t}\sum_{n \geq 0}\frac{n(4re^t)^n}{(1+re^t)^{2n}}J_t^n dt + \frac{2re^t}{(1+re^t)^2}\sum_{l \geq 0}\frac{(4re^t)^l}{(1+re^t)^{2l}}J_t^l 
\sum_{n \geq 0}  \frac{(1/2)_n}{n!}\frac{(4re^t)^n}{(1+re^t)^{2n}} \,dt   \\& + \frac{1}{2}\sum_{l \geq 0}\frac{l(4re^t)^l}{(1+re^t)^{2l}} J_t^l\sum_{n \geq 0}  \frac{(1/2)_n}{(n+1)!}\frac{(4re^t)^{n+1}}{(1+re^t)^{2n+2}} \,dt  \,
\end{align*} 

Then, using that (cf \cite{Askey}) : 
\begin{equation*} {}_1\mathscr{F}_0(a,z) : = \sum_{p=0}^{\infty}(a)_p \frac{z^p}{p!} = (1 - z)^{-a}, \, |z| <1, \end{equation*}
one has for all $t$ such that $re^t < 1$ : \begin{equation*}
\sum_{n \geq 0}  \frac{(1/2)_n}{n!}\frac{(4re^t)^n}{(1+re^t)^{2n}} \,dt  =   \left(1-\frac{4re^t}{(1+re^t)^2}\right)^{-1/2} = \frac{1+re^t }{1-re^t}
\end{equation*} and similarly, from $\int_0^u dz/\sqrt{1-z} = 2 - 2\sqrt{1-u} $,  
\begin{equation*}
\frac{1}{2}\sum_{n \geq 0}  \frac{(1/2)_n}{(n+1)!}\frac{(4re^t)^{n+1}}{(1+re^t)^{2n+2}} = 1- \left(1-\frac{4re^t}{(1+re^t)^2}\right)^{1/2} = \frac{2re^t}{1+re^t}
\end{equation*}
The result follows from an easy computation. $\hfill \blacksquare$
\begin{cor}Let $T_k$ denotes the Tchebycheff polynomial of the first kind (cf \cite{Askey}) defined by :
\begin{equation*}
T_k(x) = \cos(k\arccos(x)), \, k \geq 0, \, x \in ]-1,1[ \end{equation*}
Thus the process $S(k)$ defined by $S_t(k) := e^{kt}T_k(2J_t - P)$ is a $\mathscr{J}_t$-free martingale.
\end{cor} 
{\it Proof } : The generating function of these polynomials is given by (cf \cite{Askey}) : 
\begin{equation*}
L(x,z) := \sum_{k=0}^{\infty}T_k(x) z^k = \frac{1-zx}{1-2zx + z^2}, \quad, |z| < 1 \end{equation*}
Letting $z=re^t$ with $0 < r< e^{-t} < e^{-s}$ for $s < t$, then $L(2J_t-P, re^t)$ and $L(2J_s - P, re^s)$ converge, and using that $R$ is a $\mathscr{J}_t$-martingale, one gets : 
\begin{equation*}
\sum_{k=0}^{\infty}\tilde{\Phi}(e^{kt}T_k(2J_t - P)|\mathscr{J}_s) r^k = \tilde{\Phi}(R_t| \mathscr{J}_s) = R_s = \sum_{k=0}^{\infty}T_k(2J_s - P) e^{ks} r^k
\end{equation*} 
Taking the derivative of both sides at $r=0$, we are done.   
\section{On The Cauchy transform of the free Jacobi process}
 Under the assumptions of the previous section, one has : 
 \begin{equation*}\tilde{\Phi}(J_t) = (\tilde{\Phi}(J_0)-\theta)e^{-t} + \theta \end{equation*}
 Then, a similar computation as in Proposition (\ref{inj}) gives for all $u$ in the unit disk:
 \begin{align*}
 &\tilde{\Phi}(\log(P-uJ_t)) = -\sum_{n \geq 1}\frac{u^n}{n}\tilde{\Phi}(J_t^n) = -u \tilde{\Phi}(J_t) - \sum_{n \geq 2}\frac{u^n}{n}\tilde{\Phi}(J_t^n) \\& 
 =u(\tilde{\Phi}(J_0)-\tilde{\Phi}(J_t)) + \tilde{\Phi}(\log(P-uJ_0)) + \int_0^t \tilde{\Phi}\left(\sum_{n \geq 2}(uJ_s)^n\right) ds - \theta u \int_0^t \tilde{\Phi}\left(\sum_{n \geq 2}(uJ_s)^{n-1}\right) ds  \\& -\lambda \theta  \int_0^t\sum_{n \geq 2}\sum_{k=0}^{n-2}u^n\tilde{\Phi}(J_s^{n-k-1})\tilde{\Phi}(J_s^k(P-J_s))ds \\&
 =u(\tilde{\Phi}(J_0)-\theta)(1- e^{-t}) +  \tilde{\Phi}(\log(P-uJ_0))+ \int_0^t\tilde{\Phi}(u^2J_s^2(P-uJ_s)^{-1})ds 
 \\& - \theta u\int_0^t\tilde{\Phi}(uJ_s(P-uJ_s)^{-1})ds   -\lambda \theta u \int_0^t \tilde{\Phi}\left(\sum_{n \geq 0}(uJ_s)^{n+1}\right)\tilde{\Phi}\left(\sum_{k \geq 0}(uJ_s)^k(P-J_s)\right)ds 
 \end{align*} 
Using the fact that $u^2J_s^2 = (uJ_s - P)(uJ_s+P) + P$ and $uJ_s = (uJ_s - P) +P$, we get : \begin{align*}
&\tilde{\Phi}(\log(P-uJ_t)) = u(\tilde{\Phi}(J_0)-\theta)(1- e^{-t}) + \tilde{\Phi}(\log(P-uJ_0)) + (1-\theta u)\int_0^t \tilde{\Phi}((P-uJ_s)^{-1})ds  
\\& + \theta ut- \int_0^t \tilde{\Phi}((P+uJ_s))ds  - \lambda \theta u \int_0^t \tilde{\Phi}(uJ_s(P-uJ_s)^{-1})\tilde{\Phi}((P-J_s)(P-uJ_s)^{-1})ds 
\\&= \tilde{\Phi}(\log(P-uJ_0)) + (1 - \theta u) \int_0^t \tilde{\Phi}((P-uJ_s)^{-1})ds - t  \\& -\lambda \theta u \int_0^t [\tilde{\Phi}((P-uJ_s)^{-1}) - 1][1 + (u-1)\tilde{\Phi}(J_s(P-uJ_s)^{-1})]ds
\\& = \tilde{\Phi}(\log(P-uJ_0)) + (1 - \theta u - \lambda \theta u) \int_0^t \tilde{\Phi}((P-uJ_s)^{-1})ds - (1-\lambda \theta u)t \\& -\lambda \theta (u-1) \int_0^t \tilde{\Phi}((P-uJ_s)^{-1})\tilde{\Phi}(uJ_s(P-uJ_s)^{-1})ds  + \lambda \theta (u-1) \int_0^t \tilde{\Phi}(uJ_s(P-uJ_s)^{-1})ds 
\\&= \tilde{\Phi}(\log(P-uJ_0)) + (1 - \theta u - \lambda \theta u + 2\lambda \theta(u-1)) \int_0^t \tilde{\Phi}((P-uJ_s)^{-1})ds - (1-\lambda \theta u)t \\& -\lambda \theta(u-1)t  -\lambda \theta (u-1) \int_0^t  \tilde{\Phi}^2((P-uJ_s)^{-1})ds 
\\& = \tilde{\Phi}(\log(P-uJ_0)) + (1 - \theta u + \lambda \theta(u-2)) \int_0^t \tilde{\Phi}((P-uJ_s)^{-1})ds - (1- \lambda \theta)t  
\\&\lambda \theta (u-1) \int_0^t \tilde{\Phi}^2((P-uJ_s)^{-1})ds  \end{align*}
Setting $h_t(u) = \tilde{\Phi}((P-uJ_t)^{-1})$, then : \begin{equation*}h_t(u) = \frac{1}{u}G_t\left(\frac{1}{u}\right)\end{equation*} where $G_t$ denotes the Cauchy (or the G-) transform of the $FJP$'s law defined on \\ $\mathbb{C}^+$ by : \begin{equation*} G_t(z) = \int_0^1\frac{\mu_t(dx)}{z-x} \end{equation*} 
Furthermore, we have: \begin{equation*}-\frac{d}{du}\tilde{\Phi}(\log(P-uJ_t)) = \tilde{\Phi}(J_t(P-uJ_t)^{-1}) = \frac{1}{u}(h_t(u) - 1)\end{equation*} 
Thus, we get: \begin{align*}
h_t(u) &= h_0(u) + \theta(1-\lambda)\int_0^t u h_s(u)ds -(1 - \theta u + \lambda \theta(u-2)) \int_0^t uh'_s(u)ds \\& +\lambda \theta \int_0^t u h_s^2(u)ds + 2\lambda \theta(u-1) \int_0^t u h_s(u) h'_s(u) ds \end{align*}or equivalently : \begin{align*}
G_t(1/u) &= G_0(1/u) +\theta(1-\lambda)u\int_0^t G_s(1/u)ds +  \lambda \theta \int_0^t G_s^2(1/u)ds + (1 - \theta u + \lambda \theta(u-2))\times \\& \int_0^t\left[G_s(1/u) + \frac{1}{u}G'_s(1/u)\right]ds + 2\lambda \theta\left(\frac{1-u}{u^2}\right)\int_0^t\left[uG_s^2(1/u) + G'_s(1/u)G_s(1/u)\right] \\&
= G_0(1/u) + (1-2\lambda \theta)\int_0^t G_s(1/u)ds + \lambda \theta\left(\frac{2}{u} - 1\right)\int_0^t G_s^2(1/u)ds \\&+ \frac{1-\theta u +\lambda \theta(u-2)}{u}\int_0^tG'_s(1/u)ds 
+\frac{2\lambda \theta(1-u)}{u^2}\int_0^tG_s(1/u)G'_s(1/u)ds 
\end{align*}As a consequence, $G$ satisfies the p. d. e. : \begin{pro}
\begin{align*} G_t(z) &= G_0(z) + (1-2\lambda \theta)\int_0^t G_s(z)ds  + \lambda \theta(2z - 1)\int_0^t G_s^2(z)ds  \\&+ ((1-2\lambda \theta)z - \theta(1-\lambda))\int_0^tG'_s(z)ds + 2\lambda \theta z(z-1)\int_0^tG_s(z)G'_s(z)ds  \end{align*}\end{pro}
\section{Conclusion and open questions}
Throughout this paper, one can see that the study of the free Jacobi process  is more easier than the study of the complex matrix Jacobi process in the sense that, though both cases belong to a non commutative context, we get more powerful results on the law of the process in the infinite dimensional case, namely, the recurrence formula for the moments and the Cauchy transform. Yet, one would like to deepen our results and get much more precision on the law, that is, compute the discrete and continuous parts as done in the stationary case. 
For instance, let us interest in the Dirac mass at $1$, then $-y\Im(G(1+iy))$ involves $\Im(yG^{'}(1+iy))$ which can be unbounded. Hence, one can not intertwin limit and integral signs.

\end{document}